\documentclass[leqno]{article}
\usepackage{amssymb, amsmath, amsfonts}
\date{}

\newtheorem{theorem}{Theorem}[section]

\newtheorem{lemma}[theorem]{Lemma}
\newtheorem{proposition}[theorem]{Proposition}
\numberwithin{equation}{section}

\def\sqr#1#2{{\vcenter{\vbox{\hrule height.#2pt
    \hbox{\vrule width.#2pt height#1pt \kern#1pt
    \vrule width.#2pt}
    \hrule height.#2pt}}}}

\def\ecarre{\; \mathchoice\sqr56\sqr56\sqr{4.3}5\sqr{3.7}5}

\def\al{\alpha}
\def\De{\Delta}
\def\de{\delta}
\def\La{\Lambda}
\def\la{\lambda}
\def\eps{\varepsilon}

\def\om{\omega}

\def\pa{\partial}

\def\Dom{\hbox{Dom}}

\def\div{{\rm div}}
\def\Re{{\rm Re}\,}

\def\bu{\bar{u}}

\def\td{\tilde d}

\font\bbb=msbm10

\def\R{\hbox{{\bbb R}}}
\def\C{\hbox{{\bbb C}}}

\def\cL{{\cal L}}

\def\Dom{\, {\rm Dom}\, }

\def\Lip{\, {\rm Lip}\, }
\def\Liploc{\,{\rm Lip}_{\textit{loc}}\,}

\def\supp{\hbox{supp}\;}

\def\msni{\medskip\noindent}
\def\bsni{\bigskip\noindent}
\def\ni{\noindent}

\def\ms{\medskip}
\def\bs{\bigskip}

\begin{document}

\title{Essential self-adjointness\\
 for semi-bounded 
magnetic Schr\"odinger operators\\ on non-compact manifolds}

\author{Mikhail Shubin\footnote{Research partially supported by NSF grant DMS-9706038}}
\maketitle
\begin{abstract}
We prove essential self-adjointness for semi-bounded below magnetic
Schr\"odinger operators on complete Riemannian manifolds with
a given positive smooth measure which is fixed independently of the metric.
Some singularities of the scalar potential are allowed.

This is an extension of the Povzner--Wienholtz--Simader theorem.
The proof uses the scheme of Wienholtz but requires 
a refined invariant integration by parts technique,
as well as a use of a  family of cut-off functions
which are constructed by a non-trivial smoothing procedure due to Karcher.

\end{abstract}

\section{Introduction}\label{S:intro}

Let $(M,g)$ be a Riemannian manifold (i.e. $M$ is a $C^\infty$-manifold,
$(g_{jk})$ is a Riemannian metric on $M$), $\dim M=n$.
We will always assume for simplicity that $M$ is connected. We will also assume
that we are given a {\it positive smooth measure} $d\mu$ i.e.
a measure which has a $C^\infty$ positive density $\rho(x)$
with respect to the Lebesgue measure $dx=dx^1\dots dx^n$
in any local coordinates $x^1,\dots,x^n$, so we will write
$d\mu=\rho(x)dx$.
This measure may be completely independent of the Riemannian
metric, but may of course coincide with the canonical measure $d\mu_g$
induced by the metric (in this case $\rho=\sqrt{g}$ where
$g=\det(g_{jk})$, so locally $d\mu_g=\sqrt{g}dx$).

The main purpose of this paper is to study essential
self-adjointness of semi-bounded below magnetic Schr\"odinger operators in
$L^2(M)=L^2(M,d\mu)$.

Denote $\La^p_{(k)}(M)$ the set of all $k$-smooth (i.e. of the class $C^k$)
complex-valued $p$-forms on $M$. We will write $\La^p(M)$ instead
of $\La^p_{(\infty)}(M)$.
A {\it magnetic potential} or {\it vector potential} is a
real-valued 1-form $A\in \La^1_{(1)}(M)$. So in local coordinates
$x^1,\dots, x^n$ it can be written as
$$
A=A_jdx^j,
$$
where $A_j=A_j(x)$ are real-valued $C^1$-functions of the
local coordinates, and we use the standard Einstein summation convention.

The usual differential can be considered as a first order
differential operator
$$
d:C^\infty(M)\longrightarrow \La^1(M).
$$
We will also need a deformed differential
$$
d_A:C^\infty(M)\longrightarrow \La^1_{(1)}(M),\qquad u\mapsto du+iuA,
$$
where $i=\sqrt{-1}$.

The Riemannian metric $(g_{jk})$ and the measure $d\mu$ induce an inner
product in the spaces of smooth forms with compact support
in a standard way.
In particular, this inner product on functions has the form
$$
(u,v)=\int_M u\bar v d\mu,
$$
where the bar over $v$ means the complex conjugation.

For smooth forms
$\alpha=\alpha_j dx^j, \beta=\beta_k dx^k$ denote
$$
\langle\alpha,\beta\rangle=g^{jk}\alpha_j\beta_k,
$$
where $(g^{jk})$ is the inverse matrix to $(g_{jk})$.
So the result $\langle\alpha,\beta\rangle$ is a scalar function on $M$.
Then for $\alpha,\beta$ with compact support we have
$$
(\alpha,\beta)=\int_M \langle\alpha,\bar\beta\rangle d\mu,
$$
where
$$
\bar\beta=\bar\beta_k dx^k.
$$
Using the inner products in spaces of smooth functions and 1-forms
with compact support we can define the completions of these
spaces. They are Hilbert spaces which we will denote
$L^2(M)$  for functions and $L^2\La^1(M)$ for 1-forms.
These spaces depend on the choice of the metric $(g_{jk})$ and
the measure $d\mu$. However we will skip this dependence in the
notations of the spaces for simplicity of notations. This will not
lead to a confusion because both metric and measure will be fixed
through the whole paper unless indicated otherwise.

The corresponding local spaces will be denoted $L^2_{loc}(M)$
and $L^2_{loc}\La^1(M)$ respectively. These spaces do not depend
on the metric or measure. For example $L^2_{loc}(M)$ consists of
all functions $u:M\to\C$ such that for any local coordinates
$x^1,\dots,x^n$ defined in an open set $U\subset M$ we have $u\in L^2$
with respect to the Lebesgue measure $dx^1\dots dx^n$
on any compact subset in $U$. Similarly the space $L^p_{loc}(M)$
is well defined for any $p$ with $1\le p\le\infty$.

Formally adjoint operators to the differential operators
with sufficiently smooth coefficients are well defined
through the inner products above. In particular,
we have an operator
$$
d_A^*:\La^1_{(1)}(M)\longrightarrow C(M),
$$
defined by the identity
$$
(d_Au,\om)=(u,d_A^*\om), u\in C^\infty_c(M), \om\in\La^1_{(1)}(M).
$$
(Here $C_c^\infty(M)$ is the set of all
$C^\infty$ functions with compact support on $M$.)

Therefore we can define the magnetic Laplacian  $\De_A$
(with the potential $A$) by the formula
$$
-\De_A=d_A^* d_A: C^\infty(M)\longrightarrow C(M).
$$
Now the main object of our study will be the
{\it magnetic Schr\"odinger operator}
\begin{equation}\label{Hmag}
H=H_{A,V}=-\De_A+V,
\end{equation}
where $V\in L^2_{loc}(M)$ i.e. $V$ is a measurable
locally square-integrable function which is called {\it electric potential} or {\it scalar potential}. 
We will always assume $V$ to be real-valued. Then $H$ becomes a symmetric
operator in $L^2(M)$ if we consider it on the domain
$C_c^\infty(M)$. In this paper we will assume that
$H_{A,V}$ is semi-bounded below on $C_c^\infty(M)$
i.e. there exists a constant $C\in\R$ such that
\begin{equation}\label{E:Hsemibound}
(H_{A,V}u,u)\ge -C(u,u),\quad u\in C_c^\infty(M).
\end{equation}
We will impose the following local condition on $V$:

\medskip\noindent
$(H)$ \qquad $V=V_++V_-$ where $V_+\ge 0$, $V_-\le 0$, $V_+\in L^2_{loc}(M)$ 
and $V_-\in L^p_{loc}(M)$

\qquad\
with $p=n/2$ if $n\ge 5$, \ $p>2$ if $n=4$, and $p=2$ 
if $n\le 3$.

\msni
Our main result is the following
\begin{theorem}\label{T:self-semi} Let the manifold
$(M,g)$  be complete,
$A\in \La^1_{(1)}(M)$, $V$ satisfies the condition $(H)$ above,
and the corresponding magnetic
Schr\"odinger operator $H_{A,V}$ is semi-bounded below
on $C_c^\infty(M)$.
Then $H_{A,V}$ is essentially self-adjoint.
\end{theorem}

{\bf Remark 1.} If we assume that $V\in L^\infty_{loc}(M)$, then 
instead of $A\in C^1(M)$ it is sufficient to assume 
that $A\in \Lip_{loc}(M)$, i.e $A$ is locally Lipschitz.

\medskip
{\bf Remark 2.} The requirement on $p$ in the condition $(H)$ is almost optimal.
Indeed, we must require that $V\in L^2_{loc}(M)$ if we wish $H_{A,V}$ 
to be defined on $C_c^\infty(M)$. This is the only requirement which is imposed 
for $n\le 3$; the requirement $p>2$ in case $n=4$ is only slightly stronger.
As to the requirement $p=n/2$ in case $n\ge 5$, it can not be replaced by $p=n/2-\eps$
with $\eps>0$.  This was shown by B.~Simon even in $\R^n$ and without magnetic field
(see \cite{Simon} or \cite{Reed-Simon}, Example 4 in  Ch.X.2): the operator
$-\De-\alpha/|x|^2$ on $C_c^\infty(\R^n)$ with a real parameter $\al$  
is bounded from below if and only if
$\al\le (n-1)(n-3)/4+1/4$  and essentially self-adjoint if and only if
$\al\le (n-1)(n-3)/4-3/4$. However the requirement $V_-\in L^p_{loc}(M)$
can be replaced by weaker requirements formulated in less explicit terms,
e.g. Stummel classes \cite{Stummel} and domination requirements
(see e.g. \cite{Simader78}).

\ms
{\bf Remark 3.}
For the usual semi-bounded below Schr\"odinger operator $H=-\De+V(x)$ in $\R^3$ with a continuous 
potential $V$ the essential self-adjointness was 
conjectured by I.M.~Glazman and proved by
A.Ya.~Povzner (\cite{Povzner}, Theorem 6 in Ch.I). Independently a much more
general result (which includes in particular magnetic Schr\"odinger operators 
with sufficiently regular coefficients in $\R^n$) was obtained by E.~Wienholtz 
\cite{Wienholtz}.  The Wienholtz proof is much simpler and for the simplest case
of the Schr\"odinger operator without magnetic field
it is  also reproduced in the book of I.M.~Glazman \cite{Glazman}. 
Further  improvements for operators in $\R^n$ and in its open subsets 
are due to H.~Stetk\ae r-Hansen \cite{Stetkaer-Hansen},
J.~Walter \cite{Walter}, and C.~Simader \cite{Simader78}.  
We will use the method of Wienholtz when we treat the case of locally bounded
potentials $V$ and the method of C.~Simader \cite{Simader78} for more singular $V$.

\medskip
{\bf Example.} Let us give an example which shows that the magnetic field can
contribute to the fulfillment of the semi-boundedness condition 
(\ref{E:Hsemibound}) for $H_{A,V}$ so that the corresponding operator
$H_{0,V}$ (with the magnetic field removed) is not essentially
self-adjoint. 

In this example we will take $M=\R^2$ with the standard flat metric, so
the magnetic potential is $A=A_1dx^1+A_2dx^2$. The magnetic field is then a 2-form
$B=dA=B_{12}dx^1\wedge dx^2$. Let us write $B$ instead of $B_{12}$ for simplicity
of notation. Of course, changing the order of $x^1$ and $x^2$ would replace 
$B$ by $-B$.

Using simple uncertainty principle type arguments given e.g. in 
\cite{Colin} or \cite{Ivrii}, we can see that 
\begin{equation*}
H_{A,V}\ge B+V \qquad \text{and} \qquad H_{A,V}\ge -B+V,
\end{equation*}
where the inequalities are understood in the sense of quadratic forms.
Assume now that $V\in L^2_{loc}(\R^2)$  and either $B+V$ or $-B+V$ is semi-bounded
below.  Then due to Theorem \ref{T:self-semi} the operator $H_{A,V}$ is essentially
self-adjoint. This can happen in particular when $V\to -\infty$ fast enough so that
$H_{0,V}$ is not essentially self-adjoint, e.g. when $V(x)=-|x|^\alpha$ with 
$\alpha> 2$ (see e.g. \cite{Berezin-Shubin}, Example 1.1 in Ch.3).

\section{Algebraic preliminaries}\label{S:alg}

We will start by considering
the operator $d^*$, which is formally adjoint to $d$, so
$d^*:\La^1_{(1)}(M)\to C(M)$. This operator is related with the divergence
of vector fields. Let $v$ be a smooth vector field on $M$.
Denote by $\om_v$ the 1-form
corresponding to $v$ i.e. locally $\om_v=(\om_v)_jdx^j$ where
$$
(\om_v)_j=g_{jk}v^k.
$$
Vice versa, for any smooth 1-form $\om$ we
will denote by $v_\om$ the corresponding vector field, so locally
$v_\om=v_\om^k\pa/\pa x_k$ where
$$
v_\om^k=g^{kj}\om_j.
$$
Then we will define the divergence of $v$ by the formula
\begin{equation}\label{div-d*}
{\div}\,v=-d^*\om_v.
\end{equation}
Equivalently we can write
\begin{equation}\label{d*-div}
d^*\om=-{\div}\, v_\om.
\end{equation}
A straightforward calculation shows that in local coordinates
\begin{equation}\label{divergence}
{\div}\, v=\frac{1}{\rho}\frac{\pa}{\pa x^i}(\rho v^i),
\quad v=v^i\frac{\pa}{\pa x^i}.
\end{equation}
It follows from (\ref{div-d*}) that ${\div}\,v$
(as given by (\ref{divergence})) does not depend on the choice
of local coordinates but only on the metric and measure. On the
other hand (\ref{divergence}) implies that ${\div}\, v $ does not
depend on the metric (even though it is not immediately seen
from (\ref{div-d*})).

We have the following Leibniz rule for $d^*$
(or, equivalently, for the divergence):
\begin{equation}\label{d*-Leibniz}
d^*(f\om)=fd^*\om-\langle df,\om\rangle, \quad f\in C^1(M),
\ \om\in\La^1_{(1)}(M).
\end{equation}

For the Laplacian $\De$ (on functions) we have
\begin{equation*}
\De u= -d^*d u={\div}\,(\nabla u), \qquad u\in C^2(M),
\end{equation*}
where $\nabla u$ means the
gradient of $u$ associated with $g$, i.e. the vector field which
corresponds to $du$ and is given in local coordinates as
$$
\nabla u=g^{jk}\frac{\pa u}{\pa x^j}\frac{\pa}{\pa x^k}.
$$

Let us identify the magnetic potential $A$ with the multiplication
operator
$$
A:C^\infty(M)\longrightarrow \La^1_{(1)}(M).
$$
Then the formally adjoint operator $A^*$ is a substitution operator
of the vector field $v_A$ into 1-forms, or in other words
\begin{equation}\label{A*}
A^*\om=\langle A,\om\rangle=g^{jk}A_j\om_k.
\end{equation}
This gives us a formula for $d_A^*$:
\begin{equation}\label{dA*}
d_A^*\om=(d^*-iA^*)\om=-{\div}\,v_\om-i\langle A,\om\rangle.
\end{equation}
It follows that
\begin{equation}\label{dA*product}
d_A^*(f\om)=fd^*\om-\langle df,\om\rangle-if\langle A,\om\rangle,
\quad f\in C^1(M), \ \om\in \La^1_{(1)}(M).
\end{equation}
The following Leibniz rules for $d_A^*$
immediately follow:
\begin{eqnarray*}
&d_A^*(f\om)&=fd_A^*\om-\langle df,\om\rangle,\\
&d_A^*(f\om)&=fd^*\om-\langle d_A f,\om\rangle,
\end{eqnarray*}
where $f,\ \om$ are as in (\ref{dA*product}).

Using these formulas, we can write an explicit expression
for the magnetic
Laplacian $\De_A=-H_{A,0}$. Namely,
\begin{eqnarray*}
&-\De_Au&=d_A^*d_Au=(d^*-iA^*)(du+iAu)\\
&&=d^*du-iA^*du+id^*(Au)+A^*Au\\
&&=-\De u-i\langle A,du\rangle -i\,{\div}\,(uv_A)+\langle A,A\rangle u\\
&&=-\De u -2i\langle A, du\rangle+(id^*A+|A|^2)u.
\end{eqnarray*}
Hence we obtain the following  expression for
the magnetic Schr\"odinger operator (\ref{Hmag}):
\begin{equation}\label{HAV-explicit}
H_{A,V}u=-\De u -2i\langle A, du\rangle+(id^*A+|A|^2)u+Vu.
\end{equation}

On the other hand using the expressions (\ref{divergence}) and
(\ref{A*}) for the divergence and the operator $A^*$ we easily obtain
that in local coordinates
\begin{equation}\label{HAV-local}
H_{A,V}u=-\frac{1}{\rho}\left(\frac{\partial}{\partial x^j}+iA_j\right)
\left[\rho g^{jk}
\left(\frac{\partial}{\partial x^k}+iA_k\right)u\right]+Vu,
\end{equation}
or in slightly different notations
\begin{equation*}
H_{A,V}u=\frac{1}{\rho}(D_j+A_j)
[\rho g^{jk}(D_k+A_k)u]+Vu,
\end{equation*}
where $D_j=-i\partial/\partial x_j$.

\medskip
\textbf{Remark.} A similar operator in $\R^n$ (with $\rho=1$)
was considered
by T.~Ikebe and T.~Kato \cite{Ikebe-Kato},
K.~J\"orgens \cite{Jorgens}, M.S.P.~Eastham, W.D.~Evans and
J.B.~McLeod \cite{Eastham},
A.~Devinatz \cite{Devinatz}
in the space $L^2(\R^n, dx)$ where $dx$ is the standard Lebesgue
measure on $\R^n$. The general operator of the form (\ref{HAV-local})
on manifolds was studied by H.O.Cordes \cite{Cordes, Cordes2}.
In this generality it includes some natural geometric situations
(in particular the case $\rho=\sqrt{g}$).

\section{Preliminaries on the Lipschitz analysis and the Stokes formula
 on a Riemannian manifold}\label{S:Lip}

Let $(M,g)$ be a Riemannian manifold. A function $f:M\to \R$ is
called {\it a Lipschitz function with a Lipschitz constant $L$} if
\begin{equation}\label{Lipschitz-def}
|f(x)-f(y)|\le Ld(x,y),\quad x,y\in M,
\end{equation}
where $d(x,y)$ means the Riemannian distance between $x$ and $y$.
We will denote the space of all Lipschitz functions on $M$ by $\Lip(M)$.
This space depends on the choice of the Riemannian metric on $M$.
The space of all locally Lipschitz functions on $M$ will be denoted
$\Lip_{loc}(M)$. This space does not depend on the Riemannian metric
on $M$.

By the well known Rademacher theorem,  (\ref{Lipschitz-def}) implies that 
$f$ is differentiable almost
everywhere and
\begin{equation}\label{Lipschitz-d}
|df|\le L
\end{equation}
with the same constant $L$. Here $|df|$ means the length of the
cotangent vector $df$ in the metric associated with $g$. The
corresponding partial derivatives of the first order coincide 
with the distributional derivatives. 
Vice versa if $df\in L^\infty(M)$, for the
distributional differential $df=(\pa f/\pa x^j)dx^j$, then $f$ can
be modified on a set of measure $0$ so that it becomes a Lipschitz
function.

The estimate (\ref{Lipschitz-d}) can be also rewritten in the form
\begin{equation}\label{Lipschitz-grad}
|\nabla f|\le L
\end{equation}
(again with the same constant $L$).

In local form (in open subsets of $\R^n$) these facts are
discussed e.g. in the book of V.~Mazya \cite{Mazya}, Sect.1.1. The
correspondence between constants in (\ref{Lipschitz-def}),
(\ref{Lipschitz-d}) and (\ref{Lipschitz-grad}) is straightforward.

The Lipschitz vector fields, differential forms etc. are defined
in an obvious way.

The formulas (\ref{div-d*}), (\ref{d*-div}), (\ref{divergence}),
(\ref{d*-Leibniz}), (\ref{dA*}) apply to Lipschitz vector fields
and forms instead of smooth ones.

We will also need local Sobolev spaces $W^{m,2}_{loc}$ on $M$
for arbitrary integer $m$.
We need these spaces for functions, vector fields and differential forms. 
For simplicity let us consider functions first.
If $m\ge 0$ then the space $W^{m,2}_{loc}(M)$ consists of
functions
$u\in L^2_{loc}(M)$ such that
their derivatives of the order $\le m$ in local
coordinates also belong to $L^2_{loc}$ in these coordinates.
(The functions which coincide almost everywhere are
identified.)
Denote also by $W^{m,2}_{comp}(M)$ the space of functions which
belong to $W^{m,2}_{loc}(M)$ and have a compact support.

If $m<0$
then $W^{m,2}_{loc}(M)$ is a dual space to
$W^{-m,2}_{comp}(M)$ and it consists of all
distributions which can be  locally represented as sums of
derivatives of order $\le -m$ of functions from $L^2_{loc}$.

These definitions obviously extend to vector fields and differential forms. 

We will  need the Stokes formula, or rather the divergence formula
for Lipschitz vector fields $v$ on $M$ in the following simplest
form:

\begin{proposition}\label{Stokes}
Let $v=v(x)$ be a vector field which is in $W^{1,2}_{comp}$ on
$M$. Then $$ \int_M{\div}\, v\; d\mu=0. $$
\end{proposition}

The proof of the Proposition can be easily reduced to the case
when $v$ is supported in a domain of local coordinates. After that
we can use mollification (regularization) of $v$ to approximate
$v$ by smooth vector fields. A more advanced statement which does not require
a compact support and includes a boundary integral, 
can be proved for Lipschitz vector fields
(\cite{Mazya}, Sect. 6.2). 

\medskip
Again using mollifiers we easily see that 
the formulas (\ref{div-d*}), (\ref{d*-div}), (\ref{divergence}),
(\ref{d*-Leibniz}), (\ref{dA*}) apply to functions, vector fields
and forms from $W^{1,2}_{loc}$ instead of smooth ones.

\section{Cut-off functions}\label{S:geom}

In the proofs of the main results in the next section we will need
a sequence of compactly supported cut-off functions with Lipschitz gradients, 
such that the gradients are  uniformly small. Here we will follow H.~Karcher \cite{Karcher}
to establish the existence of such functions on any complete Riemannian manifold,
and they will be in fact $C^\infty$-functions.

\begin{proposition}\label{P:cut-off} 
Let $(M,g)$ be a complete Riemannian manifold. Then there exists
a sequence of functions
$\phi_N: M\to \R$, $N=1,2,\dots,$ with the following properties:

\msni
$(a)$ \qquad $\phi_N\in C^\infty_c(M)$;

\msni
$(b)$ \qquad $0\le \phi_N(x)\le 1, \quad x\in M$,\quad $N=1,2,\dots$;

\msni
$(c)$ \qquad for every compact $K\subset M$ there exists $N_0>0$ such that
$\phi_N=1$ on $K$

\qquad if $N\ge N_0$;

\msni
$(d)$ \qquad $\eps_N:=\underset{x\in M}{\sup}|\nabla\phi_N(x)|\to 0
\quad \textit{as} \quad N\to\infty.$
 
\end{proposition}

\ms
{\bf Proof.}
Note that for any complete Riemannian manifold $(M,g)$
it is very easy to construct a sequence of compactly supported functions
$\psi_N\in\Lip(M)$, $N=1,2,\dots,$ satisfying the conditions
$(b)$, $(c)$, $(d)$ above (if we substitute them for $\phi_N$ there). 
For example, we  can take
\begin{equation}\label{E:phi-N}
\psi_N(x)=\chi(N^{-1}d(x,x_0)),
\end{equation}
where $x_0\in M$ is fixed, 
 $\chi\in C_c^\infty(\R)$, $\chi(r)=1$ if $r\le 1/3$, $\chi(r)=0$ if $r\ge 2/3$,
and $0\le\chi(t)\le 1$ for all $r\in \R$. In this case clearly
$|\nabla\psi_N|\le C/N$.

 However it is not clear how
to satisfy (a), and it is even not  immediately clear how to make 
$\nabla\phi_N\in\Lip(M)$. But there are many manifolds where this
is easily possible, e.g. we  have $\psi_N\in C^\infty(M)$ if $M$ is $\R^n$
(with the flat metric),
the hyperbolic space, or generally  any manifold with an empty cut-locus, 
so that the function $x\mapsto d(x,x_0)$ is in $C^\infty(M)$ if $x\ne x_0$.

\ms
More generally, in the construction above we can replace
the distance function $d(x)=d(x,x_0)$ by a {\it regularized
distance}: a smooth function $\td:M\to \R$ such that $\td\ge 0$ and
$$
C^{-1}d(x)-C_1\le \td(x)\le Cd(x)+C_1
$$
with some positive constants $C,C_1$. Such a function $\td\in C^\infty(M)$
can be easily constructed on any manifold of bounded geometry
(see e.g. the construction given in \cite{Shubin}). Subtler
arguments by J.~Cheeger and M.~Gromov \cite{Cheeger-Gromov-91}
(which are based on a result of U.~Abresch \cite{Abresch} about smoothing
of Riemannian metrics,  I.~Yomdin's theorem which is a quantitative
refinement of the Sard Lemma - see \cite{Gromov} , pp. 123--124,
and some arguments from \cite{Cheeger-Gromov-85}),
allow to construct such regularized distance on any complete
Riemannian manifold with a bounded sectional curvature
(without any restrictions on the injectivity radius, which
are part of the usual definition of bounded geometry).

In the general case the result easily follows
by use of a H.~Karcher's mollifiers construction \cite{Karcher},
applied to the family $\psi_N$ from (\ref{E:phi-N}). Let us recall this construction.

Let us choose a point $m\in M$ and a small ball $B(m;\rho)$ centered at $m$ with
the radius $\rho>0$, so that this ball is geodesically convex and the exponential map
$$
\exp_m:T_mM\to M
$$
restricted to the euclidean ball $D(0;\rho)\subset T_m(M)$ is 
a diffeomorphism of $D(0;\rho)$ onto $B(m;\rho)$. We will identify $B(m;\rho)$ with 
$D(0;\rho)$ via $\exp_m^{-1}$ and construct mollifiers (or mollifying kernels)
\begin{equation}\label{E:Phi-rho}
\Phi_\rho(m,y)=\chi\left(\frac{1}{\rho}d(m,y)\right)
\left(\int_{B(m;\rho)}\chi\left(\frac{1}{\rho}d(m,x)\right)d_mx\right)^{-1},
\end{equation}
where $\chi$ is the same function as above, $d_mx$ is the euclidean volume in $B(m;\rho)$ 
(coming from $T_mM$ via the exponential map).

Choosing a compact $K\subset M$, we see that $\Phi_\rho(m,y)$ is well defined for all
$m\in K$ and arbitrary $y\in M$ provided  
$0<\rho<\rho_0=\rho_0(K)$.
Clearly $\Phi_\rho(\cdot,\cdot)\in C^\infty(U\times M)$ for a neighborhood $U$ of $K$,
$\Phi_\rho(m,y)=const$ near the diagonal $m=y$ and $\Phi_\rho(m,y)=0$ 
if $d(m,y)\ge\rho$.

H.~Karcher applied the mollifiers (\ref{E:Phi-rho}) to smooth maps $M\to \hat M$
for another Riemannian manifold $\hat M$. To this end he used the Riemannian center of mass
on $\hat M$. We will only need the case $\hat M=\R$ where the construction and 
arguments become much simpler (but still not trivial). Taking a locally integrable 
function $f:M\to\R$, we can define the mollified functions (depending on $\rho>0$) by
\begin{equation*}
f_\rho(m)=\int_M f(x)\Phi_\rho(m,x)d_mx.
\end{equation*}
Assuming for simplicity that $f$ has a compact support, $\supp f\subset K$ with a compact 
$K\subset M$, we see that $f_\rho\in C^\infty(M)$ if $\rho<\rho_0(K)$. It is also clear that
$f_\rho=0$ outside of the $\rho$-neighborhood of $K$. 

Now let us apply this to $f=\psi_N$ taking $\rho=\rho_N$ sufficiently small, and denote
the resulting mollified function by $\phi_N$, i.e. $\phi_N=(\psi_N)_{\rho_N}$. Then the 
sequence $\phi_N$, $N=1,2,\dots$, satisfies the conditions $(a)$, $(b)$ and $(c)$ above. 

It remains to see that the functions $\phi_N$ satisfy $(d)$ as well. To this end we can 
use Theorem 4.6 of Karcher \cite{Karcher}.  It implies  that if $f$ is 
a Lipschitz function with  the Lipschitz constant $L$ i.e. (\ref{Lipschitz-def}) holds,
and the sectional curvature varies in a finite interval $[\de,\De]$
in a $\rho$-neighborhood of $\supp f$, then 
$$
|f_\rho(x)-f_\rho(y)|\le L(1+L^2\cdot C(\de,\De)\rho^2)d(x,y).
$$
Hence for a Lipschitz function $f$ with a compact support, we can choose $\rho$ so small
that $f_\rho$ is Lipschitz with the Lipschitz constant $2L$. In this case we will have
$|\nabla f_\rho|\le 2L$ everywhere.  Since the Lipschitz constant of $\psi_N$ is $O(1/N)$,
the condition $(d)$ for $\phi_N$ immediately follows. $\square$

\section{Proof of Theorem \ref{T:self-semi}}\label{S:main}

In this section we will always write $H$ instead of $H_{A,V}$
for simplicity of notations.

Let $H_{min}$ and $H_{max}$ be the minimal and maximal operators
associated with the differential expression (\ref{Hmag}) for $H$
in $L^2(M)$.
Here $H_{min}$ is the closure of $H$ in $L^2(M)$ from the initial
domain $C_c^\infty(M)$, $H_{max}=H_{min}^*$ (the adjoint operator
to $H_{min}$ in $L^2(M)$). Clearly
$$
\Dom(H_{max})=\{u\in
L^2(M)|\;Hu\in L^2(M)\},
$$
where $Hu$ is understood in the sense
of distributions.

The essential self-adjointness of $H$ means that $H_{min}=H_{max}$.

For simplicity of exposition we treat the case of a locally bounded scalar
potential $V$ first.
The requirements on the vector potential $A$ can be 
slightly relaxed in this case.

\subsection{Locally bounded scalar potentials}\label{S:semi-loc-bound}

To establish the equality $H_{min}=H_{max}$  we need some information
about the domain of $H_{max}$. We will start with a simple lemma establishing
necessary local information in the simplest case $V\in L^\infty_{loc}(M)$.

\begin{lemma}\label{local-lemma}
Assume 
that $A\in \Lip_{loc}(M)$,
$V\in L^\infty_{loc}(M)$ and $u\in \Dom(H_{max})$. Then 
$u\in W^{2,2}_{loc}(M)$. 
\end{lemma}

\textbf{Proof.} We will repeat an argument given in
\cite{Berezin-Shubin}, Appendix 2, proof of Theorem 2.1.

Assume that $u\in \Dom(H_{max})$. Due to (\ref{HAV-explicit}) this
means that $u\in L^2(M)$ and
$$
-\De u -2i\langle A, du\rangle+(id^*A+|A|^2)u+Vu=f\in L^2(M),
$$
where $\De u$ and $\langle A, du\rangle$
are understood in the sense of distributions, so
a priori $\De u\in W^{-2,2}_{loc}(M)$,
$\langle A, du\rangle\in W^{-1,2}_{loc}(M)$. Note also that
$(id^*A+|A|^2)u+Vu\in L^2_{loc}(M)$. It follows from the local
elliptic regularity theorem applied to $-\De$ that
$u\in W^{1,2}_{loc}(M)$.

This already implies that $\langle A, du\rangle\in L^2_{loc}(M)$.
Applying the local elliptic regularity theorem again we see that
$u\in W^{2,2}_{loc}(M)$.
$\ecarre$

\medskip
\textbf{Remark.} Lemma \ref{local-lemma} is certainly not new,
though I had difficulty to find a statement which would exactly
imply it. More general equations are considered e.g.
by D.~Gilbarg and N.S.~Trudinger (\cite{Gilbarg-Trudinger},
Theorem 8.10), but
with  a stronger a priori requirement $u\in W^{1,2}$.

\begin{theorem}\label{T:self} Let us assume that the manifold
$(M,g)$ is complete,
$A\in \Liploc(M)$,
$V\in L^\infty_{loc}(M)$ and the corresponding magnetic
Schr\"odinger operator $H_{A,V}$ is semi-bounded below
on $C_c^\infty(M)$
i.e. (\ref{E:Hsemibound}) holds.
Then $H_{A,V}$ is essentially self-adjoint.
\end{theorem}

\textbf{Proof.} 
Note that the smoothness requirements on $A,V$ imply that the operator
$H_{A,V}$ is well defined on $C_c^\infty(M)$
and maps this space into $L^2(M)$ (see Sect.\ref{S:Lip}),
as well as on $L^2(M)$ (which it maps to the space
of distributions on $M$).

Adding $(C+1)I$ to $H_{A,V}$
we can assume that $H_{A,V}\ge I$ on $C_c^\infty(M)$
i.e.
\begin{equation*}
(H_{A,V}u,u)\ge (u,u),\quad u\in C_c^\infty(M).
\end{equation*}
If this is true, then it is well known (see e.g. \cite{Glazman})
that the essential self-adjointness of $H_{A,V}$ is equivalent
to the fact that the equation
\begin{equation*}
H_{A,V}u=0
\end{equation*}
has no non-trivial solutions in $L^2(M)$ (understood in the sense of
distributions).

Assume that $u$ is such a solution. First note that it is in
$W^{2,2}_{loc}(M)$  due to Lemma \ref{local-lemma}.

Let us take a cut-off function $\phi_N$ on $M$ from Proposition \ref{P:cut-off}.

Then denoting $u_N=\phi_N u$ we see that $u_N$ is in the domain
of the minimal operator associated with $H_{A,V}$, hence
\begin{equation}\label{E:semi-cut}
\|u_N\|^2\le (H_{A,V}u_N,u_N).
\end{equation}

Now we will prove an identity which will be useful not only in this proof
but in extending the result to singular scalar potentials.

Let us calculate $H_{A,V}(\phi u)$ for arbitrary functions $u,\phi$
such that $u\in W^{2,2}_{\textit{loc}}(M)$ and $\phi\in C^1(M)$
has a locally Lipschitz gradient. We will use the Leibniz
type formulas from Sect.\ref{S:alg}. Applying $d_A^*$ to
$$
d_A(\phi u)=\phi d_Au+ud\phi,
$$
we obtain
$$
d_A^*d_A(\phi u)=\phi d_A^*d_Au-2\langle d\phi,d_Au\rangle+ud^*d\phi,
$$
hence
\begin{equation}\label{E:HAV-product}
H_{A,V}(\phi u)=\phi H_{A,V}u-2\langle d\phi,d_Au\rangle-u\De\phi.
\end{equation}
Now let us additionally assume that $\phi$ is real-valued and has
a compact support. Then multiplying \eqref{E:HAV-product} by $\phi \bu$
and integrating over $M$ (with respect to the chosen measure $d\mu$
with a positive smooth density) we get
$$
(H_{A,V}(\phi u), \phi u)=(\phi H_{A,V}u,\phi u)-
\int_M[2\langle d\phi,du\rangle+2i\langle A,d\phi\rangle
u+u\De\phi]\phi\bu d\mu.
$$
Adding this formula with the complex conjugate one and dividing
by 2, we see that the term with $A$ under the integral sign
cancels, so using Proposition \ref{Stokes} we obtain
\begin{eqnarray*}
&(H_{A,V}(\phi u), \phi u)&=\Re (\phi H_{A,V}u,\phi u)-
\int_M[\langle \phi d\phi,\bu du+ud\bu\rangle+
|u|^2\phi\De\phi]d\mu\\
&&=\Re (\phi H_{A,V}u,\phi u)-
\int_M[\langle \phi d\phi,d(|u|^2)\rangle+
|u|^2\phi\De\phi]d\mu\\
&&=\Re (\phi H_{A,V}u,\phi u)-
\int_M\left(|u|^2d^*(\phi d\phi)+
|u|^2\phi\De\phi\right)d\mu.
\end{eqnarray*}
Since
$$
d^*(\phi d\phi)=\phi d^*d\phi-\langle d\phi,d\phi\rangle
=-\langle d\phi,d\phi\rangle-\phi\De\phi,
$$
we finally obtain the desired identity
\begin{equation}\label{E:HAV-phi-u}
(H_{A,V}(\phi u), \phi u)=\Re (\phi H_{A,V}u,\phi u)
+\int_M |d\phi|^2 |u|^2d\mu.
\end{equation}

\medskip
To use this identity in our proof assume that $H_{A,V}u=0$. This implies
$$
(H_{A,V}(\phi u), \phi u)=
\int_M |d\phi|^2 |u|^2d\mu.
$$
Now taking $\phi=\phi_N$ and applying the estimate
\eqref{E:semi-cut}, we obtain
$$
\|\phi_N u\|^2\le \int_M |\nabla\phi_N|^2 |u|^2 d\mu.
$$
In particular, for any compact $K\subset M$ we obtain
for $N\ge N_0(K)$:
$$
\int_K |u|^2d\mu\le \int_M |\nabla\phi_N|^2 |u|^2 d\mu
\le \eps_N \int_M|u|^2 d\mu.
$$
If now $u\in L^2(M,d\mu)$, then taking limit as $N\to\infty$,
we see that $u=0$ on $K$, hence $u\equiv 0$. $\ecarre$

\subsection{Singular scalar potentials}\label{SS:sing-pot}

Now we will consider magnetic Schr\"odinger 
operators $H_{A,V}$ on a complete Riemannian manifold $(M,g)$,
so that the conditions of Theorem \ref{T:self-semi} are satisfied.
In particular, we will assume that $A\in\La^1_{(1)}(M)$ but we will not require
that $V$ is locally bounded.  

\textbf{Proof of Theorem \ref{T:self-semi}.} 
1. Let us choose a relatively compact coordinate neighborhood $U$ in $M$ with
coordinates $x^1,\dots, x^n$ which are defined in a neighborhood of $\bar U$.

Let $\De_0$ denote the flat Laplacian in these coordinates. Then
due to the standard elliptic estimates  the norms 
\begin{equation*}
\|\De_0u\|+\|u\| \quad\text{and} \quad \|H_{A,0}u\|+\|u\|
\end{equation*}
are equivalent on  functions $u\in C_c^\infty(U)$. On the other hand
if we denote the bottoms of the spectra of the Friedrichs extensions of
$-\De_0$ and $H_{A,0}$ in $U$ by  $\la_0$ and $\la_A$ respectively, then 
$\la_0>0$ and also $\la_A>0$ due to the diamagnetic inequality 
(see e.g. \cite{Kato1, Simon3} or \cite{Lieb-Loss}, Sect. 7.21). 
It follows that 
\begin{equation*}
\|u\|\le\la_0^{-1}\|\De_0u\|, \quad \|u\|\le\la_A^{-1}\|H_{A,0}u\|,
\end{equation*}
for any $u\in C_c^\infty(U)$, hence there exists $C>0$ such that
\begin{equation}\label{E:HA0=De}
C^{-1}\|\De_0 u\|\le \|H_{A,0}u\|\le C\|\De_0 u\|, \quad u\in C_c^\infty(U).
\end{equation}

Now let us recall that it follows from $(H)$ that $V_-$
has  $\De_0$-bound $\eps>0$ on $C_c^\infty(U)$ for arbitrarily small $\eps$
(see Theorem X.20 and Corollary of Theorem X.21 from \cite{Reed-Simon}), i.e.
\begin{equation}\label{E:V-dom-De}
\|V_-u\|\le \eps\|\De_0u\|+C_\eps\|u\|, \quad u\in C_c^\infty(U).
\end{equation}
Using (\ref{E:HA0=De}) we see that (\ref{E:V-dom-De}) is equivalent to 
a similar estimate with $\De_0$ replaced by $H_{A,0}$:  
\begin{equation}\label{E:V-dom-HA0}
\|V_-u\|\le \eps\|H_{A,0}u\|+C_\eps\|u\|, \quad u\in C_c^\infty(U).
\end{equation}

\ms
2. We would like to extend the inequality \eqref{E:V-dom-HA0} to functions 
$u\in C_c^\infty(M)$ under the condition that 
$V_-\in L^p_{comp}(M)$ with $p$ as in $(H)$ (with $C_\eps$
depending on $V_-$). To this end we need the following estimate
\begin{equation}\label{E:dA-by-HA0}
\|d_A u\|\le \eps \|H_{A,0}u\|+C_\eps \|u\|, \quad u\in C_c^\infty(M).
\end{equation}
An equivalent form of \eqref{E:dA-by-HA0} is
\begin{equation*}
\|d_A u\|^2\le \eps \|H_{A,0}u\|^2+\tilde C_\eps \|u\|^2, \quad u\in C_c^\infty(M),
\end{equation*}
which holds because  due to the Cauchy-Schwarz inequality
\begin{equation*}
\|d_A u\|^2=(H_{A,0}u,u)\le \|H_{A,0}u\|\cdot\|u\|
\le \eps \|H_{A,0}u\|^2+\frac{1}{4\eps}\|u\|^2, \quad u\in C_c^\infty(M).
\end{equation*}
By taking closure we see that \eqref{E:dA-by-HA0} holds in fact
for all $u\in\Dom(H_{A,0})$ where the domain is understood as the domain of minimal
or maximal operators (which coincide due to \cite{Shubin3} or Theorem \ref{T:self} above).

\ms
3. Assuming that $V_-\in L^p_{comp}(M)$, let us choose functions 
$\psi_1,\dots,\psi_N\in C_c^\infty(M)$ 
such that 

\ms
(a) $\supp \psi_j\subset U_j$ for a relatively compact coordinate  neighborhood $U_j$, $j=1,\dots, N$.

\ms
(b) $\sum_{j=1}^N\psi_j=1$ in a neighborhood of $\supp V_-$.

\msni
Using  \eqref{E:V-dom-HA0}, we obtain for any $\eps>0$
\begin{equation}\label{E:Vu-by-sum}
\|V_-u\|\le\sum_{j=1}^N \|V_-(\psi_j u)\|\le \eps\sum_{j=1}^N \|H_{A,0}(\psi_j u)\|
+C_\eps \|u\|.
\end{equation}
Now we can use \eqref{E:HAV-product} to conclude that
\begin{equation}\label{E:}
\|H_{A,0}(\psi_j u)\|\le C_1(\|H_{A,0}u\|+\|d_A u\|+\|u\|), u\in C_c^\infty(M).
\end{equation}
This again holds for any $u\in \Dom(H_{A,0})$ due to the arguments given above in part 2
of this proof. We obtain now from \eqref{E:Vu-by-sum} that
\begin{equation}\label{E:V-dom-HA0-global}
\|V_-u\|\le \eps\|H_{A,0}u\|+C_\eps\|u\|, \quad u\in C_c^\infty(M),
\end{equation}  
under the condition that $V_-\in L^p_{comp}$.

\ms
4. Define 
$V_-^{(N)}(x)=V_-(x)$ on $\supp \phi_N$, and $V_-^{(N)}(x)=0$ otherwise.
(Here $\phi_N$ is the function from Proposition \ref{P:cut-off}.)
Then (\ref{E:V-dom-HA0}) holds for $V_-^{(N)}$. It follows from 
Theorem \ref{T:self} and from the Kato-Rellich
perturbation theorem (Theorem X.12 in \cite{Reed-Simon}) that the operator
$H_{A,V_-^{(N)}}=H_{A,0}+V_-^{(N)}$ is essentially self-adjoint.

Now we can  use the Kato inequality technique
(see \cite{Kato1} or \cite{Reed-Simon}, especially Theorem X.33, 
and also generalization to operators on manifolds and in sections
of vector bundles developed by H.~Hess, R.~Schrader and D.A.~Uhlenbrock
\cite{Hess77, Hess80}), and the perturbation arguments from the proofs
of Theorems  X.28, X.29 from \cite{Reed-Simon}, to prove that the operator
$H_N=H_{A,V_++V_-^{(N)}}=H_{A,0}+V_++V_-^{(N)}$ is essentially self-adjoint
for any $N=1,2,\dots$. 

Note that  the use of the Kato inequality in the last step requires that 
$A\in C^1$, rather than $\Liploc$ (see \cite{Kato1} where the Friedrichs' mollifiers 
technique \cite{Friedrichs44} is used; this technique requires the derivatives of $A$
to be continuous).

\ms
5.   In what follows we will write $H$ instead of $H_{A,V}$.
Note that for any fixed $u\in\Dom(H_{max})$
\begin{equation}\label{E:ineq}
|(u,H(\phi_N f))|=|(Hu,\phi_N f)|\le C\|f\|, \quad f\in C_c^\infty(M).
\end{equation}
Similarly to (\ref{E:HAV-product}) we have
$$
H(\phi_N f)=\phi_N Hf -2\langle d\phi_N,d_Af\rangle-f\De\phi_N,
$$
hence
$$
(\phi_N u, Hf)=2(u,\langle d\phi_N,d_Af\rangle)+(u, f\De\phi_N)+(u,H(\phi_N f)),
$$
and using (\ref{E:ineq}) we conclude that
\begin{equation*}\label{E:estimate}
|(\phi_N u,Hf)|\le C\left(\|df\|+\|f\|\right), \quad f\in C_c^\infty(M),
\end{equation*}
with the constant $C$ depending on $u$, $H$ and $\phi_N$ (but not on $f$).
Since the left hand side depends only on the restriction of $u$ to
a neighborhood of $\supp\phi_N$, we can also write
\begin{equation}\label{E:N-estimate}
|(\phi_N u,H_Nf)|\le C\left(\|df\|+\|f\|\right), \quad f\in C_c^\infty(M).
\end{equation}

\ms
6. Our next goal is to establish that $\Dom(H_{max})\subset W^{1,2}_{loc}(M)$.
It is enough to prove that (\ref{E:N-estimate}) implies
that $\phi_Nu\in W^{1,2}_{loc}(M)$. We will repeat the arguments from
\cite{Simader78}. Denote $v=\phi_Nu$, so $v\in L^2(M)$.

By the standard domination argument we have 
\begin{equation*}
|(V_-^{(N)}f,f)|\le a\|df\|^2+C\|f\|^2, \quad f\in C_c^\infty(M),
\end{equation*}
with an arbitrarily small $a>0$ and $C$ depending on $a$
or, equivalently,
\begin{equation}\label{E:dom-A-estimate}
|(V_-^{(N)}f,f)|\le a\|d_Af\|^2+C'\|f\|^2, \quad f\in C_c^\infty(M).
\end{equation}
Indeed, \eqref{E:dA-by-HA0} means 
the operator domination relation $V_-^{(N)}<< H_{A,0}$ 
which in turn implies the same domination relation for the corresponding 
quadratic forms (see Theorem X.18 in \cite{Reed-Simon}), i.e. 
\eqref{E:dom-A-estimate} with  arbitrarily small $a>0$.

Choosing an arbitrary $\lambda>0$, we obtain
\begin{equation}\label{E:H+lambda}
\begin{align}
&((H_N+\lambda)f,f)=\|d_Af\|^2+(V_+f,f)+(V_-^{(N)}f,f)+\lambda\|f\|^2\\
&\ge (1-a)\|d_Af\|^2+(\lambda-C')\|f\|^2\notag\\
&\ge (1-a)\notag\|df\|^2+ (\lambda-C'')\|f\|^2. \notag
\end{align}
\end{equation}
Now let us choose here $\lambda>C''$. Taking closure, we see that the estimate 
\eqref{E:H+lambda}  holds for all $f$ in the domain of the closure
of $H_N$ understood as the operator with the domain $C_c^\infty(M)$. 
It is a standard fact that this closure coincides with $H_N^{**}=(H_N^*)^*$.  However
since $H_N$ is essentially self-adjoint, we have $H_N^{**}=H_N^*$ and the domain 
$D_N=\Dom(H_N^{**})$  coincides with the domain of the corresponding 
maximal operator $H_N^*$, 
i.e. with the set of all  $f\in L^2(M)$ such that $H_Nf\in L^2(M)$ where $H_Nf$ 
is understood in the sense of distributions.
In particular, (\ref{E:H+lambda}) holds for all $f\in D_N$.  

Clearly, $H_N$ is semi-bounded below.  Therefore for sufficiently large $\lambda>0$ 
the operator $H_N^*+\la:D_N\to L^2(M)$ is bijective. Hence for any $\phi\in C_c^\infty(M)$
supported in the domain of some local coordinates $x^1,\dots,x^n,$ and for any 
$j\in\{1,\dots,n\}$ we can find $f_j\in D_N$ such that 
$(H_N+\lambda)f_j=\partial_j^*\phi$ where $\partial_j=\partial/\partial x^j$
and $\pa_j^*$ means the formally adjoint operator with respect to the inner product
induced by the given measure in the chosen coordinate neighborhood.
It follows that for any $\eps>0$
\begin{equation}\label{E:H+lambda-above}
|((H_N+\lambda)f_j,f_j)|=|(\pa_j^*\phi,f_j)|=|(\phi,\pa_j f_j)|
\le
\frac{\eps}{2}\|\pa_j f_j\|^2+\frac{1}{2\eps}\|\phi\|^2.
\end{equation}
Combining (\ref{E:H+lambda}) and (\ref{E:H+lambda-above}) we obtain
\begin{equation*}
\|df_j\|+\|f_j\|\le C'\|\phi\|,
\end{equation*}
with $C'$ independent of $\phi$.  Now taking $f=f_j$ in (\ref{E:N-estimate})
we obtain
\begin{equation*}
|(v,\pa_j^*\phi)|\le C''\|\phi\|.
\end{equation*}
This implies that $\pa_jv\in L^2_{loc}$ for all $j$ and $v\in W^{1,2}_{loc}$ in the  coordinate neighborhood.
Choosing a covering of $M$ by such coordinate neighborhoods we  see that
$v=\phi_Nu\in W^{1,2}_{loc}(M)$. Since $N$ was arbitrary, we see
that $u\in W^{1,2}_{loc}(M)$.
 
\medskip
7. Let us start with the identity (\ref{E:HAV-phi-u}) which was established 
in the case of a locally bounded $V$ for
all $u\in W^{2,2}_{loc}(M)$ and real-valued compactly supported $\phi$ with
a Lipschitz gradient. Let us try to relax the requirement on $u$ first,
still assuming that $V\in L^\infty_{loc}(M)$. We claim that (\ref{E:HAV-phi-u})
makes sense and holds for any $u\in W^{1,2}_{loc}(M)$. Indeed,  both sides of
(\ref{E:HAV-phi-u}) make perfect sense for any such $u$ if we understand the inner
products as dualities between $W^{-1,2}_{loc}(M)$ and $W^{1,2}_{comp}(M)$. To prove 
this identity for an arbitrary $u\in W^{1,2}_{loc}(M)$ we just need to approximate
$u$ by functions from $C_c^\infty(M)$ in the $W^{1,2}$-norm in a neighborhood of 
$\supp\phi$.

This argument works also if instead of the local boundedness of $V$ we assume 
that $V\in L^p_{loc}(M)$ where $p$ is the same as in the condition $(H)$.
Indeed, the Sobolev inequality gives a continuous imbedding of $W^{1,2}_{loc}(M)$
into $L^q_{loc}(M)$ with $q=2n/(n-2)$ if $n\ge 3$ and arbitrarily large $q<\infty$ if $n=2$.
For any $u\in W^{1,2}_{loc}(M)$ we have then $|u|^2\in L^{q/2}_{loc}(M)$ and the last
space is in a continuous dualily  with $L^p_{comp}(M)$ (by the usual
integration) due to the H\"older inequality. Therefore in this case we can
again prove the identity (\ref{E:HAV-phi-u}) for any $u\in W^{1,2}_{loc}(M)$
taking approximations by functions from $C_c^\infty(M)$.

So it remains to remove requirement $V_+\in L^p_{loc}(M)$ for $n\ge 4$ replacing it
by the inclusion $V_+\in L^2_{loc}(M)$. This can be done as follows. Let us
fix functions $u\in W^{1,2}_{loc}(M)$ and $\phi\in C^1_{comp}(M)$
with  a locally Lipschitz gradient. Then regularize $V_+$, replacing it by
$V_+^{(k)}(x)=V_+(x)$ if $V_+(x)\le k$, and $V_+^{(k)}(x)=k$ if $V_+(x)>k$;
here $k=1,2,\dots$. Then the identity (\ref{E:HAV-phi-u}) holds with $V^{(k)}=V_+^{(k)}+V_-$
instead of $V$ because $V^{(k)}\in L^p_{loc}(M)$. 
But now we can take limit as $k\to\infty$. The only terms depending on $k$ in 
(\ref{E:HAV-phi-u}) will be two identical terms 
$$
\int_M V_+^{(k)}|\phi u|^2 d\mu
$$
in the left and right hand sides.  This integral obviously has a limit 
(possibly $+\infty$) because the integrand 
converges monotonically.  By the Beppo Levi theorem this limit equals
\begin{equation}\label{E:V+-integral}
\int_M V_+|\phi u|^2 d\mu,
\end{equation}
so taking $k\to\infty$ we see that (\ref{E:HAV-phi-u}) holds for $V$.

If we only require that $u\in W^{1,2}_{loc}(M)$, then both sides of (\ref{E:HAV-phi-u})
can possibly be $+\infty$. If we know however that $u\in\Dom(H_{max})$ then the right hand
side is finite (which in fact just means the finiteness of the integral 
(\ref{E:V+-integral})). Then the left hand side is finite too.

8. Using the identity (\ref{E:HAV-phi-u}) which is now established for all
$u\in \Dom(H_{max})$, we can finish the proof of  Theorem \ref{T:self-semi}
by repeating the arguments of the proof of Theorem \ref{T:self} which follow
after this identity. $\square$

\section{Examples and further comments}\label{S:ex}

In this section we will provide several examples, further results
and some relevant bibliographical comments (by necessity incomplete).

\ms
\textbf{1.} Let us comment about the \textit{gauge invariance} for the
magnetic Schr\"odinger operators. It is easy to see that if we
replace $A$ by $A'=A+d\phi$ with a real-valued $\phi\in C^1(M)$,
such that $\nabla\phi\in\Liploc(M)$, then we have
\begin{equation*}
H_{A',V}=e^{-i\phi}H_{A,V}e^{i\phi},
\end{equation*}
both for minimal and maximal operators defined by the expression
$H_{A,V}$. Therefore it is clear that being essentially
self-adjoint is a gauge invariant property, i.e. it does not
change under any gauge transformation $A\mapsto A+d\phi$.
This well known observation was extended by H.~Leinfelder
\cite{Leinfelder} to a very general class of operators and
gauge transformations with minimal regularity conditions.
He considered the case $M=\R^n$ (with the standard metric)
but his arguments can be easily extended
to the case of arbitrary Riemannian manifolds, so we will
formulate the result for the general case.
Let us consider a class $\cL_2(M)$ which consists of operators
$H_{A,V}$ on a Riemannian manifold $(M,g)$
with $A\in L^4_{\textit{loc}}(M)$, $d^*A\in L^2_{\textit{loc}}(M)$
and $V\in L^2_{\textit{loc}}(M)$. Assume further that we have
two operators $H_{A,V}, H_{A',V}\in \cL_2(M)$ and $A'=A+d\phi$
where $\phi$ is a distribution on $M$. Then the essential self-adjointness
properties for $A$ and $A'$ are equivalent.

\ms
If $M$ has vanishing cohomology $H^1(M,\R)$ (e.g. if $M$ is simply-connected)
then  the gauge invariance above means that the essential
self-adjointness depends in fact on the magnetic field $B=dA$
(which is a 2-form or a de Rham current of degree 2) and not on
the magnetic potential $A$ itself.

\bs
\textbf{2.} Let us give some particular cases of Theorems \ref{T:self} 
 and \ref{T:self-semi}.

\begin{theorem}\label{Gaffney}
Let $(M,g)$ be a complete Riemannian manifold. Then the magnetic
Laplacian  $-\De_A=d_A^*d_A$
is essentially self-adjoint in $L^2(M,d\mu)$ for any
magnetic potential $A\in\Liploc(M)$ and any positive smooth
measure $d\mu$.
\end{theorem}

\textbf{Proof.}
Take $V\equiv 0$ and use Theorem \ref{T:self}.
$\ecarre$

Theorem \ref{Gaffney} generalizes the classical theorem by
M.~Gaffney \cite{Gaffney} which corresponds to the case when $A=0$
and $d\mu=d\mu_g$ is the Riemannian measure.

Note however that in fact the proof of Theorem \ref{T:self-semi} uses
some elements of the Gaffney proof.

N.N.~Ural'ceva \cite{Uralceva} and S.A.~Laptev
\cite{Laptev} provided examples of elliptic operators in $L^2(\R^n,dx)$
of the form
$$
\frac{\pa}{\pa x^j}\left(g^{jk}(x)\frac{\pa}{\pa x^k}\right)
$$
(with smooth positive definite matrices $(g^{jk})$) which are
not essentially self-adjoint due to the fact that the coefficients
$g^{jk}$ are ``rapidly growing". In these examples the inverse
matrix $(g_{jk})$ is vice versa ``rapidly decaying", which implies
that $\R^n$ with the metric $(g_{jk})$ is not complete.

\begin{theorem}\label{Carleman}
Let $(M,g)$ be a complete Riemannian manifold with a positive
smooth measure $d\mu$, $A\in\La^1_{(1)}(M)$,
$V\in L^2_{loc}(M)$, and $V(x)\ge -C, \ x\in M,$
with a constant $C$. Then the magnetic Schr\"odinger operator $H=-\De_A+V(x)$
is essentially self-adjoint.
\end{theorem}

In case when $M=\R^n$ (with the standard metric and measure), $A=0$
and $V\in L^\infty_{loc}(\R^n)$ this result was established independently by
T.~Carleman~\cite{Carleman} and K.~Friedrichs \cite{Friedrichs}.
(The Carleman proof is also reproduced in
the book of I.M.~Glazman \cite{Glazman}, Theorem 34 in Sect.3.)
The fact that in this case the requirement $V\in L^\infty_{loc}(\R^n)$ can be
replaced by $V\in L^2_{loc}(\R^n)$,
was established by T.~Kato \cite{Kato1} (see also \cite{Reed-Simon},
Sect. X.4). 

The work by T.~Kato was partially motivated by the paper of
B.~Simon \cite{Simon} who proved the essential self-adjointness
under an additional restriction compared with \cite{Kato1}.
The reader may consult Chapters X.4, X.5 in
M.~Reed and B.~Simon \cite{Reed-Simon}
for more references, motivations and a review.

Though the completeness requirement looks  natural in case
of semi-bounded operators, sometimes it can be relaxed
and incompleteness may be compensated by a specific behavior
of the potential (see e.g. A.G.~Brusentsev \cite{Brusentsev}
and also the references there).

We will mention a few more references which might be useful for the reader.
Reviews of different aspects of self-adjointness can be found e.g. in
\cite{Kalf, Kalf-75, Reed-Simon, Shubin2, Shubin3}.
Papers by M.~Braverman \cite{Braverman}, and M.~Lesch \cite{Lesch} contain
conditions for essential self-adjointness of some operators on sections of vector bundles.
In particular, operators considered in \cite{Lesch} generalize magnetic
Schr\"odinger operators.
Semi-bounded operators of higher order were studied by
A.G.~Brusentsev \cite{Brusentsev-85}. A.~Iwatsuka \cite{Iwatsuka} gave explicit 
conditions on the potentials of the magnetic Schr\"odinger operator in $\R^n$
(including interaction of electric and magnetic fields) which are sufficient
for the essential self-adjointness.  Different aspects of essential self-adjointness in
domains in $\R^n$ and manifolds with boundary where behavior of 
the coefficients near the boundary is relevant, were studied e.g. 
by A.G.~Brusentsev \cite{Brusentsev, Brusentsev-98}, K.~J\"orgens \cite{Jorgens},
R.~Mazzeo and R.~McOwen \cite{Mazzeo-McOwen}.
Finite speed propagation is an alternative method to prove essential 
self-adjointness (P.~Chernoff \cite{Chernoff}, A.A.~Chumak \cite{Chumak}).
I.~Oleinik discovered a new method which makes the relation between classical
and quantum completeness (the later means essential self-adjointness)
more explicit -- see \cite{Oleinik1, Oleinik2, Oleinik3,
Shubin2, Shubin3, Braverman, Lesch}.
H.~Leinfelder and C.~Simader \cite{Leinfelder-Simader} 
(see also \cite{Cycon-Froese-Kirsch-Simon}) proved the essential 
self-adjointness for the magnetic Schr\"odinger operators $H_{A,V}$ in $\R^n$
with $V\ge 0$ and with the minimal local regularity requirements on $A,V$.

\medskip
About other conditions of essential self-adjointness
for $H_{0,V}$ and $H_{A,V}$ formulated in terms of the potentials and sometimes allowing
operators which are not semi-bounded below see e.g. 
\cite{Berezanski, Faris-Lavine, Gimadislamov, Hellwig, Ismagilov, Kalf-Rofe-Beketov, 
Kato2, Levitan, Nelson, Rofe-Beketov1, Rofe-Beketov2, Rofe-Beketov3, 
Schechter, Schechter2, Simader74, Simader92, Simon2, Titchmarsh, Unell}
and references there.

\bsni
\textbf{Acknowledgments.}

\ni
\textit{
I am very grateful to U.~Abresch, M.~Braverman, J.~Br\"uning, H.~Kalf,  J.~Lott, 
O.~Milatovich, Yu.~Netrusov and F.S.~Rofe-Beketov for useful discussions and
bibliographical remarks.}

\bigskip\noindent
Department of Mathematics\\
Northeastern University\\
Boston, MA 02115, USA

\msni
E-mail: shubin@neu.edu

\bsni
AMS subject classification:

\msni
Primary 35P05, 58G25; Secondary 47B25, 81Q10

\end{document}